\documentclass[11pt, oneside]{article}   	% use "amsart" instead of "article" for AMSLaTeX format
\usepackage{geometry}                		% See geometry.pdf to learn the layout options. There are lots.
\usepackage{graphicx}				% Use pdf, png, jpg, or eps§ with pdflatex; use eps in DVI mode
\usepackage{amssymb}

\def\qed{\ \hfill$\bullet$}
\def\cqed{\ \hfill$\circ$}
\def\C{{\cal C}}
\def\E{{\cal E}}
\def\N{{\rm I\!N}}
\def\S{{\cal S}}
\def\mod{{\rm mod\ }}

\title{H-Kernels by walks}
\author{Hortensia Galeana-S\'anchez\\ hgalena@matem.unam.mx\\ Instituto de Matem\'aticas, UNAM \and Hugo Rincon-Galeana\\ hugo.rincon.galeana@gmail.com\\ Technical University of Vienna \and Ricardo Strausz\\ dino@math.unam.mx\\ Insituto de Matem\'aticas, UNAM}
\begin{document}
\maketitle

\abstract{
Let $D=(V,E)$ and $H=(U,F)$ be digraphs and consider a colouring of the arcs of $D$ with the vertices of $H$; we say that $D$ is $H$ coloured. We study a natural generalisation of the notion of kernel, as introduced by V. Neumann and Morgenstern~(1944), to prove that {\it If every cycle of $D$ is an $H$-cycle, then $D$ has an $H$-kernel by walks\/}. As a consequence of this, we are able to give several sufficient condition for the existence of $H$-kernels by walk; in particular, we solve partially a conjecture by Bai et al. in this context \cite{Bai} .
}

%\section{}
%\subsection{}

\section{Introduction}

Let $D=(V,E)$ and $H=(U,F)$ be digraphs; $D$ would be assumed not to have loops, however, $H$ may have them. We say that {\sl $D$ is $H$-coloured\/} whenever we have a function $\varsigma\colon E\to U$; that is, we coulour the arcs of $D$ with the vertices of $H$. Given such a colouring, a walk $W=z_0\dots z_k$ in $D$ is said to be an {\sl $H$-walk\/} if the sequence of colours of such a walk form a walk in $H$; that is, $\varsigma(z_0z_1)\varsigma(z_1z_2)\dots\varsigma(z_{k-1}z_k)$ is a walk in $H$. Analogously, we define {\sl $H$-trail, $H$-path\/} and {\sl $H$-cycle\/} --- unless otherwise said, every cycle considered in the sequel is a directed cycle. 

Given a digraph $D=(V,E)$, $H$-coloured, a subset of its vertices $K\subseteq V$ is said to be {\sl $H$-independent (by walks)\/} if there are no $H$-walks between any two of its elements. $K$ is said to be {\sl $H$-absorbent (by walks)\/} if for every other vertex $v\in V\setminus K$ there is an $H$-walk from it to an element of $K$. If $K$ is both, $H$-independent and $H$-absorbent, we say that it is {\sl an $H$-kernel of $D$\/}. 

%(aca el contexto)

The classical concept of {\sl kernel\/} in digraphs was introduced by John Von Neumann and Oskar Morgensten in the context of game theory~\cite{NM} where a kernel was originally named ``a~solution''; they were searching for mathematical principles applied to ``decision making''. Later on, Vassek Sv\'atal proved in~\cite{VS} that to decide the existence of a kernel is NP-complete; even in very spacial cases (e.g. planar with restricted degrees), as proved by A.S.~Fraenkel, the problem remains NP-complete~\cite{F}. One of the deepest results in the area is the celebrated theorem of Richardson~\cite{R} which claims that {\it every digraph without odd directed cycles has a kernel.\/} There has been several generalisation and variations of Richardson's theorem (e.g., see~\cite{GM} and the references therein) one of which is very important due to its applications, namely the theorem of Pierre Duchet~\cite{PD}: 

$\bullet$ {\it If every directed cycle has a symmetric arc, then the digraph has a kernel.\/}

The concept of {\sl kernel by monochromatic paths\/} ---where $H$ consists of loops only--- was introduced by Sans, Sawer and Woodrow in \cite{SSW}; since then several research has been done searching for sufficient conditions, typically in the colouring, for such a kernel to exists. Classes of much study are tournaments, bipartite, cuasi-transitive, to mention just a few... 

Later on, Sans, in a join work with Linek~\cite{SL}, consider colouring the edges of a tournament with the elements of a partial order, and introduced the concept of {\sl monote\/} paths, and the corresponding generalisation of kernel. In the same paper, they start the study of $H$-colourings replacing the partial order by a digraph. Linek continue this study, in a join work with Arpin \cite{AL}, and introduced implicitly the notion of {\sl panchromatic patern\/}~\cite{HS} --- it remains as an interesting open problem to characterize such patterns.

The {\it $H$-closure of $D$,\/} $\C_H(D)=(V(D),A)$, is the digraph which consists of the same vertices than $D$ and there is an arc $uv\in A$ iff there is an $H$-walk from $u$ to $v$. It is well known that:

$\bullet$ {\it $D$ has an $H$-kernel if and only if $\C_H(D)$ has a kernel.\/}

Let $D=(V,E_1\sqcup E_2)$ be a digraph endowed with a partition of its edges. $S\subseteq V$ is an {\it $H$-semikernel $\mod E_1$\/} if $S$ is $H$-independent and for every $z\in V\setminus S$, if there exists an $Sz$ $H$-walk, with its edges are contained in $E_2$, then there exists a $zS$ $H$-walk. 

Given another digraph $G$, we say that $D$ is {\it $G$-rainbow free\/} if $D$ does not contains a sub digraph isomorphic to $G$ with all its edges receiving different colours. In particular, we will deal with $\{C_3,P_3\}$-rainbow free digraphs; that is, they do not contain a copy of the directed cycle $C_3$, of order 3, nor a copy of the path $P_3$, of length 3, with its 3 edges of different colour.

The main contribution of this paper is:

$\bullet$ {\it If every cycle of $D$ is an $H$-cycle, then $D$ has an $H$-kernel by walks.\/}

This solves partially a conjecture by Bai et al. \cite{Bai} working with walks instead of paths --- their $H$ consists of the complete loop-less digraph. In particular, if $H$ consist of only two vertices, joined with the two possible edges, it is the same to work with walks than to work with paths. More generally, whenever the existence of and $H$-walk implies the existence of an $H$-path, our result solves their conjecture.

\section{Main Lemma}

{\bf Lemma.} {\sl Let $D$ be an $H$-coloured digraph. If every (directed) cycle in $D$ is an $H$ cycle, then every cycle in $\C_H(D)$ has a symmetric arc.}

\smallskip
{\bf Proof.} Let us suppose that, for the contrary, there is a directed cycle $\gamma = (v_0v_1\dots v_{n-1})$ in $\C_H(D)$ without symmetric arcs; further more, suppose such a cycle is of minimum length. Let $C_i=c^i_1c^i_2\dots c^i_{n_i}$ and $C_{i+1}=c^{i+1}_1c^{i+1}_2\dots c^{i+1}_{n_{i+1}}$ be $H$-walks in $D$ of minimum length which witness the existence of the arrows $v_iv_{i+1}$ and $v_{i+1}v_{i+2}$ in $\gamma$; that is, $c^i_1=v_i$, $c^i_{n_i}=v_{i+1}=c^{i+1}_1$ and $c^{i+1}_{n_{i+1}}=v_{i+2}$ ---here and in the sequel, the indices in $\gamma$ are taken mod~$n$. 

\medskip
{\it Claim 1.\/} {\sl For all $i=0,1,\dots,n-1$, the walks $C_i$ and $C_{i+1}$ do intersect only in $v_{i+1}$; i.e., $C_i\cap C_{i+1}=\{v_{i+1}\}$.}

For the contrary, let $v=c^i_j=c^{i+1}_k$ be the first vertex in $C_{i+1}$ which belongs to $C_i$ too. Observe that $(vc^i_{j+1}c^i_{j+2}\dots v_{i+1}c^{i+1}_{2}c^{i+1}_3\dots c^{i+1}_k)$ is a closed directed walk in $D$ which must contain a directed cycle $\gamma_{i+1}$ which passes through $v_{i+1}$. Since every cycle of $D$ is an $H$-cycle, $\gamma_{i+1}$ is. By the minimality of $C_i$ (analogously $C_{i+1}$), there is a unique arc $c^i_1c^i_2$ from $v_i$ in $C_i$ and a unique arc $c^i_{{n_i}-1}c^i_{n_i}$ to $v_{i+1}$ in the same walk; thus, the last arc of $C_i$ and the first of $C_{i+1}$ belong to $\gamma_{i+1}$ and form an $H$-walk. Therefore the concatenation of $C_i$ and $C_{i+1}$ is an $H$-walk in $D$. By the minimality of $\gamma$, the cycle $(v_0v_1\dots v_iv_{i+2}\dots v_{n-1})$ in $\C_H(D)$ must contain a symmetric arrow, and this has to be $v_{i+2}v_i$; let $C'_{i+2}$ be an $H$-walk minimal witness of such an arrow. A similar argument shows that the concatenation of $C_{i+1}$ and $C'_{i+2}$ is an $H$-walk from $v_{i+1}$ to $v_i$ which contradicts that $v_iv_{i+1}$ was asymmetric. This proves the claim. \cqed

\medskip
{\it Claim 2.\/} {\sl For all $i\not=j$, if $j\not=i\pm1$, then the walks $C_i$ and $C_{j}$ are non-intesecting.}

For the contrary, let $\E:=\{C_i \mid \exists j\not=i : C_i\cap C_j\not=\emptyset\}$ and suppose it is not empty. Let $f\colon\E\to\N$ be defines as
	$$f(C_i)=\min\{j\not=i\mid C_i\cap C_j\not=\emptyset\}.$$
Observe that, by the previous claim, for all $C_i\in\E$ it follows that $|i-f(C_i)|\geq2$.

Let $C_i$ be such that $k:=|i-f(C_i)|$ is as small as possible; w.l.o.g., we may suppose that $f(C_i)=i+k$. Let $W:=\{v_i,v_{i+1},\dots,v_{i+k}\}$ denote those vertices in $\gamma$ corresponding to the $H$-walks $W':=\{C_i,C_{i+1},\dots,C_{i+k}\}$. Observe that the only walks in $W'$ that intersect are $C_i$ and $C_{i+k}$.

Let $v=c^i_r=c^{i+k}_s\in C_i\cap C_{i+k}$, and consider the following closed walk:
	$$W'':=(c^i_r c^i_{r+1} \dots v_{i+1} c^{i+1}_1 \dots v_{i+2} \dots v_{i+k} c^{i+k}_1 \dots c^{i+k}_s).$$
Observe that the closed walk $D$ contains a directed cycle $\gamma_{i+1}$ which contains the vertex $v_{i+1}$. Since each $C_\ell$ is minimal in length, the last arc of $C_i$ and the first arc of $C_{i+1}$ must be present in $\gamma_{i+1}$. Since every cycle is an $H$-cycle, those two arcs form an $H$-path. Therefore, the concatenation of $C_i$ with $C_{i+1}$ is an $H$-walk from $v_i$ to $v_{i+2}$, and the pair $v_iv_{i+2}$ is an arc of $\C_H(D)$. Let 
	$$\gamma'=(v_0v_1\dots v_iv_{i+2}\dots v_{n-1})$$
be a cycle in $\C_H(D)$. Since $\gamma'$ is smaller than $\gamma$, it cannot be asymmetric, and the only arc which can be symmetric is $v_iv_{i+2}$. Let $C_t$ be a minimal $H$-walk representing the arc $v_{i+2}v_i$. Consider now the closed walk $C_iC_{i+1}C_t$. Such a walk contains a cycle going through $v_{i+2}$, which by hypothesis is an $H$-cycle, and the minimality of $C_{i+1}$ and $C_t$ implies that the last arc of $C_{i+1}$ and the first of $C_t$ is and $H$-path. This implies that there is an $H$-walk from $v_{i+1}$ to $v_i$ contradicting the fact that the arc $v_iv_{i+1}$ was asymmetric in $\C_H(D)$. This proves the claim.
 \cqed

Putting together both claims, we can observe that the closed walk 
	$$\hat\gamma=(C_0C_1\dots C_{n-1})$$
contains cycles going through each $v_i$ which uses the last arc of $C_{i-1}$ and the first arc of $C_i$ implying that $\hat\gamma$ is an closed $H$-walk. Therefore, the vertices $v_0 v_1 \dots v_{n-1}$ induces a complete digraph in $\C_H(D)$; in particular, the cycle $\gamma$ contains a symmetric arc, contradicting the original hypothesis and concluding the proof of the lemma.
\qed

\section{Applications.}

\bigskip
{\bf Theorem 1.} {\sl If $D$ is an $H$-coloured digraph such that every (directed) cycle is an $H$-cycle, then $D$ has an $H$-kernel by walks.}

\smallskip
{\bf Proof.} Let $\C_H(D)$ be the $H$-closure of $D$. Since $D$ has an $H$-kernel by walks if and only if $\C_H(D)$ has a kernel, it is enough to prove that $\C_H(D)$ has a kernel. Due to the Berge-Douchet theorem, it is enough to prove that every cycle in $\C_H(D)$ has a symmetric arc; but this is the content of the main lemma, which concludes the proof. \qed

\bigskip
\noindent
{\bf Theorem 2.} {\sl Let $D=(V,E)$ be an $H$-coloured digraph, let $E=E_1\sqcup E_2$ be a partition of its arcs and let $D_i=(V,E_i)$, for $i=1,2$. If every cycle contained in $D_2$ is an $H$-cycle, then $D$ has a non-empty $H$-semikernel, modulo $E_1$, which consists of a single vertex.}

\smallskip
{\bf Proof.} In seek of a contradiction, let us suppose that for every $v\in V$, $\{v\}$ is not an $H$-semikernel, modulo $E_1$. Let $v_0\in V$ be any vertex of $D$. Since $\{v_0\}$ is not an $H$-semikernel, mod $E_1$, there exists a $v_1\in V$ and a $v_0v_1$ $H$-walk contained in $E_2$ such that there is no $v_1v_0$ $H$-walk in $D$; that is, $v_1$ is the witness of $v_0$ not being an $H$-semikernel mod $E_1$. Analogously, there is a witness of $v_1$, which we denote by $v_2$, and recursively we construct $v_{i+1}$ as the witness of $v_i$. Since $D$ is finite, the sequence $v_0v_1v_2\dots v_n$ must have a vertex $v_i$ which is witness of a vertex $v_{i+k}$ and we have a closed $H$-walk $\gamma$ in $E_2$ consisting of the concatenation of their respective walks. Observe that $\gamma$ induces an asymmetric cycle in $\C_H(D_2)$. But, by hypothesis we have that every cycle of $D_2$ is an $H$-cycle, which allows us to use the lemma, contradicting the asymmetry of such a cycle. \qed

\bigskip
\noindent
{\bf Theorem 3.} {\sl Let $\S=\S(D)$ be the $H$-semikernel digraph of $D=(V,E)$ mod $E_1\subset E$. If every cycle of $D$ contained in $E_1$ is an $H$-cycle, then $\S$ is acyclic.}

\smallskip
{\bf Proof.} In seek of a contradiction, suppose that $\S$ has a cycle $\gamma=\{S_1,\dots,S_n\}$, $n>1$. First of all, observe that if for all $i\in[n]$, and each vertex $v\in S_i$ we have that $v\in S_{i+1}$, then $\gamma=\{S_1\}$, contradicting that $n>1$; therefore, we may suppose that there exists $i\in[n]$ and a vertex $v_1\in S_i$ such that $v\not\in S_{i+1}$. Since there is an arc from $S_i$ to $S_{i+1}$, then there exists a witness $v_2\in S_{i+1}$ such that there exists a $v_1v_2$ $H$-walk in $E_1$ and no $H$-walk from $v_2$ to any vertex of $S_i$. Now, observe that there must exist an $j=i+k$ minimum with the property that $v_2\in S_j$ and $v_2\not\in S_{j+1}$, since otherwise if $v_2$ is in all $S_j$ it would be in $S_i$ which implies the existence of an $H$-walk from $S_i$ to $S_i$ contradicting the $H$-independence of $S_i$. Therefore, there exists a witness $v_3\in S_{j+1}$ of the arc $S_jS_{j+1}$. Recursively, we construct the sequence $v_1v_2\dots v_m$, $m\leq n$, which is induces an asymmetric cycle in $\C_H(D_1)$, where $D_1:=(V,E_1)$. Since all cycles in $D_1$ are $H$-cycles, by the Lemma all cycles in $\C_H(D_1)$ have a symmetric arc, which is a contradiction and concludes the proof of the theorem.\qed

\bigskip
\noindent
{\bf Theorem 4.} {\sl Let $D=(V,E)$ be an $H$-coloured digraph, let $E=E_1\sqcup E_2$ be a partition of its arcs and let $D_i=(V,E_i)$, for $i=1,2$. Let us suppose that every cycle in $D_i$, for $i=1,2$, is an $H$-cycle. Furthermore, suppose that $D$ is $\{C_3,P_3\}$-rainbow free. If every $H$-walk of $D$ is contained in a $D_i$, then $D$ has an $H$-kernel by walks.}

\smallskip
The main idea behind the following proof is to choose a vertex $K$ of $\S(D)$ with exdegree zero; such a vertex exists since, by Theorem~3, $\S(D)$ is acyclic. Then, we prove that $K\subset V$ is the desired kernel.

\smallskip
{\bf Proof.} Let $K$ be a vertex of exdegree zero in $\S=\S(D)$; in search of a contradiction, suppose that $K\subset V$ is not an $H$-kernel of $D$. Since $K$ is a vertex of $\S$, then $K$ is $H$-independent (in $D$). Let $X$ be the set of vertices $v\in V$ such that there is no $vK$ $H$-walk in $D$; observe that we are supposing that such an $X$ is not empty. Since $D[X]$ is a non empty subdigraph of $D$, it follows that $D[X]$ satisfies the hypothesis of Theorem~2, and therefore there exists a vertex $x_0$ in $X$ such that $\{x_0\}$ is an $H$-semikernel, mod $E_1\downharpoonright X$, of $D[X]$. Let 
	$$T:=\{z\in K\subset V : \not\exists\ zx_0\ H{\rm-walk\ in}\ D_1\}.$$

%%Hence, for all $z\in K\setminus T$ there is a $zx_0$ H-walk in $D_1$. 

We prove now that $T\cup\{x_0\}$ is a $H$-semikernel mod $E_1$. For, since $T\subset K$, it follows that $T$ is $H$-independent by walks; and clearly $\{x_0\}$ is too. Since $x_0\in X$, there are no $x_0K$ $H$-walks, then there are no $x_0T$ $H$-walks; therefore there are no $Tx_0$ $H$-walks in $D_2$. It follows that $T\cup\{x_0\}$ is $H$-independent. We need now to show the conditional $H$-absorbence. For, let $z\not\in T\cup\{x_0\}$. Let us suppose that there is a $Tz$ $H$-walk in $D_2$ and denote it by $\alpha_1=u,\dots,a,z$. Since $K$ is a semikernel, there is an $H$-walk $\alpha_2=z,b,\dots,c,w$ with $w\in K$; we may suppose that $w\in K\setminus T$. Let $\alpha_3=w,d,\dots,x_0$ be an $H$-walk in $D_1$. The independence of $K$ implies that $azb$ is not an $H$-walk, and we may suppose that $cwd$ is not either; but, $uzwx_0$ is a $P_3$-rainbow, which contradicts the hypothesis. So, we may suppose there is a $x_0z$ $H$-walk in $D_2$; denote it by $\alpha_1=x_0,\dots,a,z$ and suppose that $z\not\in X$. It follows that there is an $H$-walk $\alpha_2=z,b,\dots,c,w$, with $w\in K$. We may suppose that $x\not\in T$, and then there is the $H$-walk $\alpha_3=w,d,\dots,x_0$. Since $x_0\in X$, then $azb$ is not an $H$-walk, and we may suppose that $bwc$ is not either; but, $x_0zwx_0$ is a $C_3$-rainbow, which contradicts again our hypothesis. Therefore, $T\cup\{x_0\}$ is a $H$-semikernel mod $E_1$ and a vertex of $\S$.

Finally, observe that there is an arc in $\S$ from $K$ to $T\cup\{x_0\}$ since $T\subset T\cup\{x_0\}$ and for each $w\in K\setminus T$ there is a $wx_0$ $H$-walk in $D_1$ and there is no $x_0K$ $H$-walk in $D$. This contradicts our choice of $K$, which is of exdegree zero, and ends the proof. \qed

The previous argument can be generalised as follows:

\bigskip
\noindent
{\bf Chorolary 5.} {\sl Let $D=(V,E)$ be an $H$-coloured digraph, let $E=A_1\sqcup A_2\sqcup\dots\sqcup A_n$ be a partition of its arcs and let $D_i=(V,A_i)$, for $i=1,2,\dots,n$. Let us suppose that every cycle in $D_i$, for $i=1,2,\dots,n$, is an $H$-cycle. Furthermore, suppose that $D$ is $\{C_3,P_3\}$-rainbow free. If every $H$-walk of $D$ is contained in a $D_i$, and ${\cal C}(D)$ is bipartite, then $D$ has an $H$-kernel by walks.}

\smallskip
{\bf Proof.} For, let $E=E_1\sqcup E_2$ be the bipartition of $E$ induced by the bipartition of ${\cal C}(D)$ and apply the previous Theorem~4.\qed

\bigskip
\noindent
{\bf Chorolary 6.} {\sl Let $D=(V,E)$ be an $H$-coloured digraph, let $E=A_1\sqcup A_2\sqcup\dots\sqcup A_n$ be a partition of its arcs and let $D_i=(V,A_i)$, for $i=1,2,\dots,n$. Let us suppose that every cycle in $D_i$, for $i=1,2,\dots,n$, is an $H$-cycle. Furthermore, suppose that $D$ is $\{C_3,P_3\}$-rainbow free. If every $H$-walk of $D$ is contained in a $D_i$, and ${\cal C}(D)$ is strongly connected without odd directed cycles, then $D$ has an $H$-kernel by walks.}

\smallskip
{\bf Proof.} For, recall that an strongly connected digraph without odd directed cycles is bipartite, and apply the previous Chorolary~5.\qed

\end{document}